\documentclass[12pt]{amsart}
\usepackage{amssymb}

\usepackage{graphicx}

\usepackage[]{amsfonts}

\def\underset#1#2{\mathrel{\mathop{\kern0pt #2}\limits_{#1}}}

\def\overset#1#2{\mathrel{\mathop{\kern0pt #2}\limits^{#1}}}

\def\couleur(#1 #2 #3)
	{
\special{ps::[end]}}

\def\sqr#1#2{{\vcenter{\vbox{\hrule height.#2pt
             \hbox{\vrule width.#2pt height#1pt \kern#1pt
             \vrule width.#2pt}
             \hrule height.#2pt}}}}

\def\st{\mathinner{\mkern1mu\raise1pt\hbox{.}				
		   \mkern1mu\raise4pt\hbox{.}
		   \mkern1mu\raise1pt\hbox{.}
		 }
         }

\def\bx#1{\setbox1=\hbox{\kern3pt{#1}\kern3pt}				
 \dimen1=\ht1 \advance\dimen1 by 3pt \dimen2=\dp1 \advance\dimen2 by 3pt
 \setbox1=\hbox{\vrule height\dimen1 depth\dimen2\box1\vrule}%
 \setbox1=\vbox{\hrule\box1\hrule}%
 \advance\dimen1 by .4pt \ht1=\dimen1
 \advance\dimen2 by .4pt \dp1=\dimen2 \box1\relax}

\def\k#1{\kern#1em}
\def\vci{\vrule  width.02em height1.47ex depth-.0ex}				
\def\11{{\rm\k{.2}\vci\k{-.37}1}}

\begin{document}
\title{A mathematical Model for Rogue Waves,\\using Saint-Venant Equations 
with Friction}
\author{Alain-Yves LeRoux , Marie-No{\"e}lle LeRoux}
\address{UNIVERSITE BORDEAUX1, Institut Math{\'e}matiques de Bordeaux, 
UMR 5251,351,Cours de la Lib{\'e}ration, 33405, Talence Cedex}
\email{Alain-Yves.Leroux@math.u-bordeaux1.fr\\}
\maketitle
\begin{abstract} {
We propose to contruct a temporary wave on the surface of the ocean, 
as a particular solution of the Saint-Venant equations with a source 
term involving the friction, whose shape is expected to mimic a rogue 
wave\ \par
}\end{abstract}
\ \par
{\hskip 1.2em}The phenomenon of Rogue Waves, or Freak Waves, is a transitory 
phenomenon which appears in the open ocean, under the shape of a gigantic 
and devastator wave, because it accumulates a significant quantity of 
energy. Many articles are devoted to them in ~\cite{ERWBrest2004}. We 
propose here an hydraulic model, making use of the Saint-Venant equations. 
The model is limited to only one dimension of space, in the direction 
of the propagation, since the side movements near the wave are uniform. 
The validity of this model requires a large wavelength, particularly 
in deep ocean. The wave itself has a rather short wavelength, but the 
whole phenomenon involves two waves of large wavelength, and relatively 
low amplitude.\ \par
{\hskip 1.8em}We denote by $H$ the mean depth of the ocean, and by $c_{s}$ 
the sonic velocity of waves inside the water (we have $c_{s}=1647\ ms^{-1}$). 
The wavelength $\lambda $ of the Saint-Venant waves must satisfy a condition 
of the form
\begin{equation} 
\lambda \ \geq \ 2\ N\ \frac{\displaystyle H\ \displaystyle \sqrt{\displaystyle 
gH}}{\displaystyle c_{s}}\ \ ,\label{RogueWave0}
\end{equation}  where $g$ is the gravity constant and $N$ the number 
of sonic interactions (back-and-forth) between the bottom and the surface 
of the ocean. This condition means that along a horizontal distance of 
a wavelength $\lambda $, there are at least $N$ such sonic interactions 
. The use of the Saint-Venant model is as more appropriate as $N$ is 
great. We usually require $N\geq 25,$ which implies for example a wavelength 
greater than  $21400\ m$ for an ocean depth of $3700m$. We suppose here, 
for simplicity, that the propagation goes from West ($x<0$) to East ($x>0$).\ 
\par
\section{The initial configuration of the wave{\hskip 1.8em}}
\setcounter{equation}{0}We denote by $q=q(x,t)$ the ocean depth at a 
point $x$ at a time $t$, and $m=m(x,t)$ the corresponding flux. We have 
$m=q\ u$, where $u$ is the water velocity. The relative velocity of the 
waves of the ocean surface is given by 
\begin{displaymath} 
c\ =\ \sqrt{gq}\ \ \ \ (=c(q))\ .\end{displaymath}  We denote by $k>0$ 
the friction coefficient of Strickler and we consider as a simplification 
purpose that the bottom is flat. In this configuration, the Saint-Venant 
model reads 
\begin{equation} 
q_{t}\ +\ m_{x}\ =\ 0\ ,
\end{equation}  
\begin{equation} 
m_{t}\ +\ 2u\ m_{x}\ +\ \left({c^{2}-u^{2}}\right) \ q_{x}\ +\ k\ \left\vert{u}\right\vert 
u\ =\ 0\ \ .
\end{equation}        We shall use the representation of waves proposed 
in ~\cite{ESource1} \author{when a source term (here $k\left\vert{u}\right\vert u$)\\ 
}
is present. The different states of the same wave are described by a 
segment of a straight line 
\begin{equation} 
m\ =\ A\ q\ -\ B\ ,
\end{equation}  in the phase plane, that is the plane $(q,m)$ here. The 
parameter $A$ is a constant which corresponds to the wave velocity, of 
the form $A=u_{ref}-c(q_{ref})$ or $A=u_{ref}+c(q_{ref})$ depending if 
the wave is travelling eastwards (sign $+$) or westwards (sign $-$)  
and for a given reference state $M_{ref}=(q_{ref},m_{ref}),$ with $m_{ref}=q_{ref}u_{ref},$ 
of course. The parameter $B$ is also a constant and is determined by 
the reference state since $B=m_{ref}+A\ q_{ref}.$ We denote by $M_{0}$ 
the state of the ocean far on the west side, and by $M_{*}$ the state 
of the ocean far on the east side. In both cases, the distances are supposed 
to be larger than the reference wavelength $\lambda $ proposed in~(\ref{RogueWave0}), 
which allows to use the Saint-Venant model. We suppose that these states 
also corresponds to reference velocity equal to zero; this assumption 
is linked to the hypothesis of a flat bottom of the ocean. That way $M_{0}$ 
represents a state $M_{0}=\left({q_{0},0}\right) $ and $M_{*}$ represents 
a state $M_{*}=\left({q_{*},0}\right) $. We suppose
\begin{displaymath} 
q_{0}\ >\ q_{*}\ \ .\end{displaymath}  A difference $q_{0}-q_{*}$ of 
a few decimeters is enough even for a wide depth of the ocean. The profile 
of the initial state is made of two branches which meet in a state $P=\left({q_{P},m_{P}}\right) $ 
for example at $x=0,$ which corresponds to the choice of the origin. 
From the state $P$ to the state $M_{*}$, that is on the East side, the 
profile is decreasing and referred to the state $M_{*}$. The corresponding 
states are so situated on the straight line 
\begin{displaymath} 
m\ =\ c_{*}\ \left({q-q_{*}}\right) \ \ ,\ \end{displaymath} with $c_{*}=c(q_{*})=\sqrt{gq_{*}}\ .$ 
As a matter of fact, we have  $A=c_{*}$ and $B=c_{*}q_{*}$ for this part 
of the wave profile. The explicit formulation of the profile is obtained, 
following ~\cite{ESW0}, by inverting the profile relation 
\begin{equation} 
\psi _{E}(q)\ =\ x\ ,\label{RogueWave1}
\end{equation}  where the index $E$ stands for East. The inverse profile 
$\psi _{E}(q)$ is determined by integrating 
\begin{displaymath} 
\psi _{E}'(q)\ =\ \frac{B^{2}-c^{2}q^{2}}{k\left\vert{Aq-B}\right\vert 
\left({Aq-B}\right) }\ \ ,\ \ \psi _{E}(q_{P})=0\ \ .\end{displaymath} 
 We need to have $q_{*}\ <\ q_{0}\ <\ q_{P}$ , in order to ensure an 
increasing profile on the West side, then a decreasing one on the East 
side. We get that way a positive flux: $Aq-B=c_{*}\left({q-q_{*}}\right) \ >\ 0$, 
and by writing  $\xi =\frac{q}{q_{*}}$ we get 
\begin{displaymath} 
\psi _{E}'(q)\ =\ \frac{1-\xi ^{3}}{k\left({1-\xi }\right) ^{2}}\ =\ \frac{1+\xi 
+\xi ^{2}}{k\left({1-\xi }\right) }\ =\ \frac{1}{k}\ \left({\ \frac{3}{1-\xi 
}\ -2\ -\xi }\right) \ ,\end{displaymath}  thus
\begin{displaymath} 
\psi _{E}(q)\ =\ -\frac{3q_{*}}{k}\ ln\left({\frac{q-q_{*}}{q_{P}-q_{*}}}\right) 
\ +\ \frac{2}{k}\ \left({q_{P}-q}\right) \ +\ \frac{1}{2kq_{*}}\ \left({q_{P}^{2}-q_{*}^{2}}\right) 
\ .\end{displaymath} By inverting this function $\psi _{E}(q)$ and using 
~(\ref{RogueWave1}) we get  $q$ as a decreasing function of $x$ which 
is equal to $q_{P}$ when $x=0\ .$      For the West side of the profile, 
with the index $W$, the reference state $M_{ref}=\left({q_{ref},m_{ref}}\right) \ $must 
correspond to a depth satisfying
\begin{displaymath} 
q_{ref}\ \geq \ q_{P}\ (\geq q_{0})\ \ ,\end{displaymath} in order to 
ensure an increasing profile. The reference velocity associated with 
the state $M_{ref}$ that is 
\begin{displaymath} 
A_{ref}\ =\ \frac{m_{ref}}{q_{ref}}\ +\ c_{ref}\ \ \ ,\ \ \ with\ \ c_{ref}\ 
=\ c(q_{ref})\ =\ \sqrt{gq_{ref}}\ \ ,\end{displaymath} is also the velocity 
of the West side profile of the wave and is positive. The straight line 
representing this West profile in the phase plane has the form
\begin{displaymath} 
m\ =\ m_{ref}\ +\ A_{ref}\ \left({q-q_{ref}}\right) \ ,\end{displaymath} 
and passes by the state $M_{0}\ =(q_{0},0)$. Hence 
\begin{displaymath} 
m_{ref}\ =\ A_{ref}\ \left({q_{ref}-q_{0}}\right) \ \ \ ,\ \ B_{ref}\ 
=\ c_{ref}\ q_{ref}\ \ .\end{displaymath}   The invert profile is described 
by a   function $\psi _{W}$ satisfying
\begin{displaymath} 
\psi _{W}'(q)\ =\ \frac{B_{ref}^{2}-c_{ref}^{2}q_{ref}^{2}}{k\ \left({A_{ref}q-B_{ref}}\right) 
^{2}}\ \ ,\ \ \psi _{W}(q_{P})\ =\ 0\ \ ,\end{displaymath} with $q_{0}\ \leq \ q_{P}\ \leq \ q_{ref}\ .$ 
The initial West side profile of the wave is then obtained by inverting, 
for any $x<0$, 
\begin{equation} 
\psi _{W}(q)\ =\ x\ .\label{RogueWave3}
\end{equation}       We set 
\begin{displaymath} 
\xi \ =\ \frac{q}{q_{ref}}\ \ \ ,\ \ \ \xi _{0}\ =\ \frac{q_{0}}{q_{ref}}\ 
\ \ ,\ \ F_{ref}\ =\ \frac{m_{ref}}{c_{ref}\ q_{ref}}\ ,\ (\ Froude\ 
number\ )\ ,\ \end{displaymath} to obtain 
\begin{displaymath} 
\psi _{W}(q)\ =\ \frac{q_{ref}}{k\ \left({F_{ref}+1}\right) ^{2}}\ \int_{\frac{q_{P}}{q_{ref}}}^{\frac{q}{q_{ref}}}{\left({-\xi 
-2\xi _{0}+\frac{\xi _{0}\left({\xi _{0}-4}\right) }{\xi -\xi _{0}}+\frac{1-\xi 
_{0}^{3}}{\left({\xi -\xi _{0}}\right) ^{2}}}\right) \ d\xi }\ \ ,\end{displaymath} 
that is
\begin{displaymath} 
\psi _{W}(q)=K\left[{\frac{q_{P}^{2}-q^{2}}{2\ q_{ref}}+2\frac{q_{0}}{q_{ref}}\left({q_{P}-q}\right) 
+\frac{q_{0}}{q_{ref}}\left({q_{0}-4q_{ref}}\right) ln\left({\frac{q-q_{0}}{q_{P}-q_{0}}}\right) 
+\frac{\left({q_{ref}^{3}-q_{0}^{3}}\right) \left({q-q_{P}}\right) }{q_{ref}\left({q-q_{0}}\right) 
\left({q_{P}-q_{0}}\right) }}\right] \end{displaymath} with 
\begin{displaymath} 
K\ =\ \frac{1}{k\ \left({F_{ref}+1}\right) ^{2}}\ \ .\end{displaymath} 
 The whole initial profile is now given by ~(\ref{RogueWave3}) for $x<0$ 
and by ~(\ref{RogueWave1}) for $x>0$.  It corresponds to an increasing 
function $q(x,0)$ for $x<0$ and a decreasing one for  $x>0$, which is 
continuous at $x=0$ where its value is $q=q_{P}$.\ \par
\section{The propagation of the wave}
\setcounter{equation}{0}\ \par
{\hskip 1.8em}The wave profile is expected to propagate eastwards, with 
the respective velocities  $A_{ref}$ and $c_{*}$ which are different 
for each part West or East. Since the west profile moves lightly faster, 
the crest of the profile will move up, at the junction of the two parts. 
The left side of this crest corresponds to the continuity of the West 
profile, extropolated for depth values $q$ going from $q_{P}$ to the 
maximal value $q_{ref}$.\ \par

\begin{figure}[h]
\begin{center}
\rotatebox{0}{\resizebox{10cm}{!}{\includegraphics{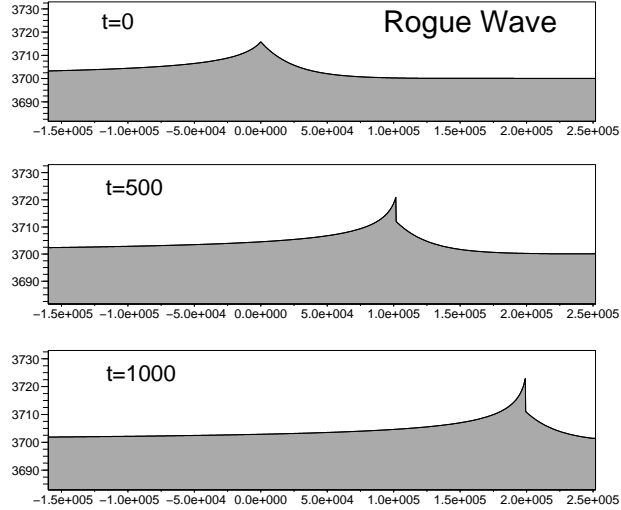}}}
\caption{Propagation of a Rogue Wave}
\label{RW1}
\end{center}
\end{figure}
\ \par
{\hskip 1.8em}The right part of the crest corresponds to a discontinuity, 
that is a shock wave, whose location is imposed by the mass conservation. 
At any time $t$ the water contained in the bump under the crest comes 
from the column of water of length $\left({A_{ref}-c_{*}}\right) t\ $shifted 
since the initial time. We denote by $q_{l}(t)$ and $q_{r}(t)$ the water 
depths on the left side (index $l$) and on the right side (index $r$) 
of the discontinuity, and by $x_{0}(t)$ its position. We always have
\begin{displaymath} 
q_{*}\ \leq \ q_{r}(t)\ \leq \ q_{P}\ \leq \ q_{l}(t)\ \leq \ q_{ref}\ 
.\end{displaymath}  As the West profile is moving with the constant velocity 
$A_{ref}$, we get it by simply inverting for any time $t$, the relation
\begin{displaymath} 
\psi _{W}(q)\ =\ x\ -\ A_{ref}t\ \ \end{displaymath} for $x<x_{0}(t)$, 
with the function $\psi _{W}$ defined above. By the same way, the East 
profile is obtained by inverting the relation
\begin{displaymath} 
\psi _{E}(q)\ =\ x\ -\ c_{*}t\ \end{displaymath} for $x>x_{0}(t)$, with 
the function $\psi _{E}$ defined above. For a given time $t$, the depths 
$q_{l}t)$ and $q_{r}(t)$, and the shock position $x_{0}(t)$ are linked 
by three conditions, entailing three equations. The first equation says 
that $q_{l}(t)$ is the value of the West profile at $x=x_{0}(t)$, that 
is 
\begin{displaymath} 
\psi _{W}(q_{l}(t))\ =\ x_{0}(t)\ -\ A_{ref\ }t\ .\end{displaymath}  
The second equation says that   $q_{r}(t)$ is the value of the East profile 
at $x=x_{0}(t)$, that is 
\begin{displaymath} 
\psi _{E}(q_{r}(t))\ =\ x_{0}(t)\ -\ c_{*}\ t\ .\end{displaymath}  The 
third equation is given by the compatibility relation of Rankine-Hugoniot
\begin{equation} 
\displaystyle \left({\displaystyle q_{l}(t)-q_{r}(t)}\right) \displaystyle 
\sqrt{\displaystyle g\frac{\displaystyle q_{r}(t)+q_{l}(t)}{\displaystyle 
2q_{r}(t)q_{l}(t)}}\ +\ \frac{\displaystyle A_{\displaystyle ref}\ q_{0}}{\displaystyle 
q_{l}(t)}\ -\ \frac{\displaystyle c_{*}\ q_{*}}{\displaystyle q_{r}(t)}\ 
=\ A_{\displaystyle ref}\ -\ c_{*}\ ,\label{EnglishGwagouFall6}
\end{equation}  which ensures the mass conservation. A dichotomy method 
running on the parameter $x_{0}(t)$ allows the simultaneous determination 
of these three parameters. For practical purposes it is however easier 
to determine $x_{0}(t)$ by checking directly the mass conservation. Let 
us consider two points $x_{1}$ and $x_{2}$ such that 
\begin{displaymath} 
x_{1}\ +\ A_{ref}\ t\ <\ x_{0}(t)\ <\ x_{2}\ +\ c_{*}\ t\ .\ \end{displaymath} 
If it is not the case in a first choice, one can increase $x_{2}$ or 
decrease $x_{1}$ sufficiently. We denote by $M_{0}$ the mass of water 
laying between the two points $x_{1}$ and $x_{2}$ at the initial time: 

\begin{displaymath} 
M_{0}\ =\ \int_{x_{1}}^{x_{2}}{q(x,0)\ dx}\ =\ \int_{x_{1}}^{0}{q_{W}(x,0)\ 
dx}\ +\ \int_{0}^{x_{2}}{q_{E}(x,0)\ dx}\ ,\end{displaymath}  where $q_{W}$ 
and $q_{E}$ are the respectives depths of the West and East profiles. 
The same mass $M_{0}$ has to be found at any time $t$ between the two 
trajectories of equations 
\begin{equation} 
x_{1}'(t)\ =\ A_{\displaystyle ref}\ -\ \frac{\displaystyle q_{\displaystyle 
ref}c_{\displaystyle ref}}{\displaystyle q_{W}\displaystyle \left({\displaystyle 
x_{1}(t)-A_{\displaystyle ref}t}\right) }\ \ ,\ x_{1}(0)\ =\ x_{1}\ ,\label{RogueWave4}
\end{equation}  and 
\begin{equation} 
x_{2}'(t)\ =\ c_{*}\ \displaystyle \left({\displaystyle 1-\frac{\displaystyle 
q_{*}}{\displaystyle q_{E}\displaystyle \left({\displaystyle x_{2}(t)-c_{*}t}\right) 
}}\right) \ \ ,\ \ x_{2}(0)\ =\ x_{2}\ ,\label{RogueWave5}
\end{equation}  which means that the relation $M(t)\ =\ M_{0}$ always 
occurs, with 
\begin{displaymath} 
M(t)\ =\ \int_{x_{1}(t)}^{x_{2}(t)}{q(x,t)\ dx}\ .\end{displaymath}  
We remark that the mass $M(t)$ may be explicitely calculated from the 
relation
\begin{displaymath} 
M(t)=\int_{q_{1}(t)}^{q_{l}(t)}{q\ \psi _{W}'(q)\ dq}\ +\ \int_{q_{r}(t)}^{q_{2}(t)}{q\ 
\psi _{E}'(q)\ dq}\ ,\end{displaymath}  where $q_{1}(t)=q_{W}(x_{1}(t))\ ,\ q_{2}(t)=q_{E}(x_{2}(t))\ .$ 
This calculation involves only primitives of rationnal fractions. More 
simply, using the computing files giving $q_{W}$ on the West side of 
a point $x_{0}$ and $q_{E}$ on the East side of this point $x_{0}$, one 
can compute 
\begin{displaymath} 
F(x_{0})\ =\ \int_{x_{1}(t)}^{x_{0}}{q_{W}\ dx\ }+\ \int_{x_{0}}^{x_{2}(t)}{q_{E}\ 
dx}\ ,\end{displaymath} and next determine $x_{0}=x_{0}(t)$ such that 

\begin{displaymath} 
F(x_{0})\ =\ M_{0}\end{displaymath}  by noticing that $F$ is an increasing 
function.\ \par

\begin{figure}[h]
\begin{center}
\rotatebox{0}{\resizebox{10cm}{!}{\includegraphics{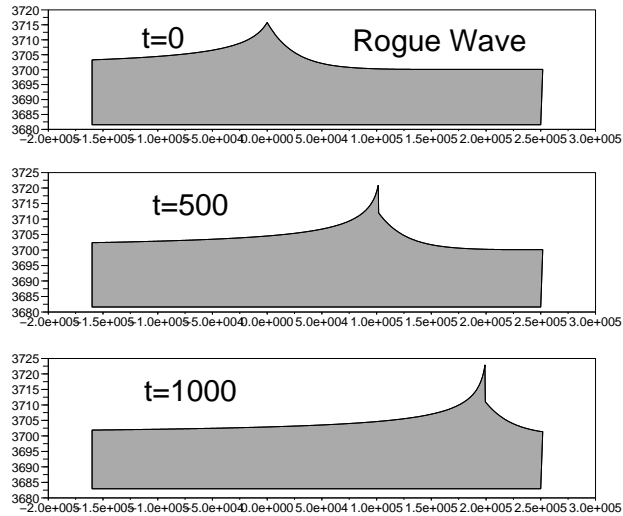}}}
\caption{Propagation of a bigger wave}
\label{RW2}
\end{center}
\end{figure}
\ \par
\ \par
{\hskip 1.8em}This last process has been used to compute the results 
shown on the figure above. For this example, we have taken $q_{*}=3700\ m$ 
and $q_{0}=3700.2\ m$, then $q_{ref}$ at its maximal value, that is $q_{ref}=3731.6737\ m$ 
here. This choice leads to the value of $q_{P}$, which is $q_{P}=3715.8087\ m$ 
here. The friction coefficient was taken equal to $k=0.45\ .$ We have 
to notice that the value of the friction coefficient is strongly linked 
to the wavelength of the profiles since the smallest frictions give the 
larger wavelengths. For terrestrial hydraulic flows (rivers or estuaries 
for example) the values of the friction coefficient is far smaller, of 
the order of $10^{-3}$ but corresponds to far less deep flows. It seems 
natural to consider that for larger depths this parameter has to be upgraded. 
At the time $t=1000$, the shock amplitude $q_{l}-q_{r}$ reaches the value 
$11\ m$, between the two cells adjascent to the computed position of 
$x_{0}(t)$ . The relative error on the mass is of order of $10^{-4}$ 
(mainly due to the meshsize), and the trajectories ~(\ref{RogueWave4}) 
and ~(\ref{RogueWave5}) were approached by a simple trapezoid formula, 
since the depth has a small variation along the trajectories $x_{1}(t)$ 
and $x_{2}(t)$ when they are chosen relatively far from the shock.\ \par
{\hskip 2.1em}By upgrading lightly the value of $q_{0}$, one get more 
important shock amplitudes. For the next example, with $q_{0}=3700.8\ m$, 
we obtain the values $q_{ref}=3763.8773,$ and $q_{P}=3731.8248$. The 
computed shock amplitude $q_{l}-q_{r}$ is equal to $33.787\ m$. The wave 
crest culminates at more than 50 meters above the sea level.\ \par
\section{Graphic interpretation{\hskip 1.8em}}
\setcounter{equation}{0}We propose a graphic interpretation in the phase 
plane $(q,m)$ using the picture below.\ \par

\begin{figure}[h]
\begin{center}
\rotatebox{0}{\resizebox{10cm}{!}{\includegraphics{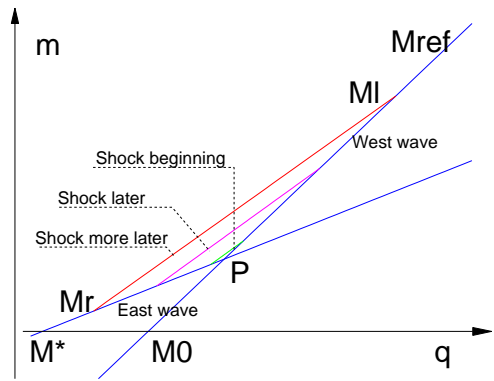}}}
\caption{The representation in the phase plane}
\label{Phase}
\end{center}
\end{figure}
\ \par
{\hskip 1.8em}The states $M_{*}$ and $M_{0}$ are represented on the $q-$axis. 
The line passing through $M_{*}$ corresponds to the East wave, of equation 
$m=c_{*}(q-q_{*})$. The line passing through $M_{0}$ corresponds to the 
West wave, of equation $m=A_{ref}(q-q_{0})$. These two lines meet at 
the point $P$. The shock wave is represented by the segment $M_{l}M_{r}$, 
with $M_{l}$ sliding along the West wave line, from $P$ to $M_{ref}$ 
and $M_{r}$ sliding along the East wave line, from $P$ to $M_{*}$. This 
segment $M_{l}M_{r}$ is not exactly a straight one, and its equation 
is obtained from the Rankine Hugoniot condition ~(\ref{EnglishGwagouFall6}). 
The maximal amplitude is reached when the first event between $M_{l}$ 
reaching $M_{ref}$ or $M_{r}$ reaching $M_{*}$ occurs. After that time 
the wave is expected to collapse. The scaling of the picture has been 
strongly modified for readability, since in the real example all those 
lines are so close to one another that it is impossible to perceive any 
difference.\ \par
\ \par
\section{Conclusion{\hskip 1.8em}}
\setcounter{equation}{0}This study shows that the outbreak of these waves 
is due to a differential in the pressure field, resulting from two waves 
with different profiles, with close but different velocities and sufficiently 
large wavelengths. The water is pushed up from this pressure effect together 
with a friction effect. The later behaviour is not studied here: one 
expects that, when $q_{l}(t)$ reaches the la value $q_{ref}$ (which occurs 
in a finite time), a backwards wave will appear, with a negative velocity, 
and will provoke strong perturbations on the West profile near the wave 
crest, which will collapse soon. However this backwards wave will have 
a shorter wavelength, probably too small to fit up with the use of the 
 Saint-Venant model, as described by the condition ~(\ref{RogueWave0}). 
\ \par
{\hskip 1.8em}Another remaining work is to fit out the different parameters 
in order to be in accordance with real world observations.  We also outline 
the strong sensitivity of the difference $q_{0}-q_{*}$ on the hight of 
the wave crest, and then on the shock amplitude. A little more important 
difference between $q_{0}$ and $q_{*}$ causes noticeably more important 
shock amplitudes. We have only proposed empirical choices of the parameters, 
in order to get realistic results showing that this study may be a suitable 
way to understand  the shaping of Rogue Waves. We emphasize the friction 
plays here a fundamental role since it allows a linear behaviuor of the 
West and the East waves.\ \par
{\hskip 1.8em}Another idea to retain is a new example of the application 
of the notion of Source Waves after the Roll Waves in channels and rivers, 
the Tidal Bore waves in the estuaries, the surf waves on the shore near 
the beach, the hurricanes and tsunamis (see ~\cite{ECras04}, ~\cite{ESW0}). 
\ \par
{\hskip 1.8em}The authors thank Michel Olagnon from Ifremer-Brest for 
some useful answers by e-mail. Reading for example  ~\cite{ENorwayWave} 
or ~\cite{EDysthe1} is also very instructive for the description of the 
phenomenon of Rogue Waves and starting some bibliography research. It 
appears that a discussion on either the linear character or the non linear 
character of the waves is spreading. In this study we put together both 
characters since linearity is introduced through the source term, that 
is the friction term, with different parameters on the two sides of the 
wave, and the nonlinearity effect is present in the shockwave, linking 
the two sides of the wave. \ \par
\ \par

\bibliographystyle{c:/TeXLive/texmf/bibtex/bst/base/plain}

\end{document}